\documentclass{gen-p-l}
\usepackage{amssymb}

\newcommand{\FF}{{\mathbb{F}}}

\newcommand{\fA}{{\mathfrak{A}}}

\newcommand{\Aut}{{\operatorname{Aut}}}
\newcommand{\Out}{{\operatorname{Out}}}

\newcommand{\SL}{{\operatorname{SL}}}
\newcommand{\PSL}{{\operatorname{L}}}
\newcommand{\OO}{{\operatorname{O}}}
\newcommand{\SU}{{\operatorname{SU}}}
\newcommand{\PSU}{{\operatorname{U}}}
\newcommand{\Atlas}{{$\mathbb{ATLAS}$}}
\newcommand{\GAP}{{\sf GAP}}
\newcommand{\Magma}{{\sc Magma}}

\newcommand{\tw}[1]{{}^#1\!}

\newtheorem{theorem}{Theorem}[section]

\theoremstyle{definition}

\theoremstyle{remark}
\newtheorem{remark}[theorem]{Remark}

\numberwithin{equation}{section}

\begin{document}

\title{Reliability and reproducibility of Atlas information}


\author{Thomas Breuer}
\address{Lehrstuhl D f{\"u}r Mathematik, RWTH Aachen, 52062 Aachen, Germany.}
\email{thomas.breuer@math.rwth-aachen.de}
\author{Gunter Malle}
\address{FB Mathematik, TU Kaiserslautern, Postfach 3049,
         67653 Kaisers\-lautern, Germany.}
\email{malle@mathematik.uni-kl.de}
\author{E.A.\ O'Brien}
\address{Department of Mathematics, University of Auckland,
Private Bag 92019, Auckland, New Zealand}
\email{e.obrien@auckland.ac.nz}


\date{\today}

\begin{abstract}
We discuss the reliability and reproducibility of much of the information
contained in the \Atlas{} of Finite Groups.
\end{abstract}

\maketitle


\section{Introduction}   \label{sec:intro}

The \Atlas{} of Finite Groups \cite{Atl}, published in 1985, contains a wealth
of information on the sporadic simple groups, their covering groups and
automorphism groups, as well as on numerous other finite simple groups
of small order. It has become an indispensable tool for researchers not
only in finite group theory but in many other areas where finite groups
play a role. In a recent letter, Jean-Pierre Serre stated that he ``can't
think of any other book published in the last 50 years which had such an
impact'', while Benedict Gross is cited as saying that
if ever the university library caught fire and one could save just one
book, the obvious choice would be the \Atlas.
In view of this, the question of reliability and reproducibility
of the results stated there is of considerable importance, particularly since
the \Atlas{} does not contain proofs of the information
it records, although for the sporadic groups it gives a list of
references from which some of the stated results are taken.
\par

In the intervening thirty years, some misprints and errors have
been found in the \Atlas. Corrections and improvements known at the time of
publication of the Atlas of Brauer Characters \cite{JLPW95} are listed in
an appendix \cite{BN95} to that book; the website \cite{AtlCorr} reports those
found later. No corrections have been added since 2002.

Few of these concern the actual character tables; apart from the
misprints, only three cases are known in which the printed tables contain
mathematical mistakes not arising simply from inconsistencies concerning
irrationalities and power maps; these concern the nearly simple groups
$2.\PSL_4(3).2_3$, $\PSU_3(11).2$ and $2.\PSU_6(2).2$;
see the more detailed comments in Section~\ref{sec:chartabs-mistakes}.
\par

Our purpose is to provide references for some of the cited results; to report
on the independent reconstruction of most of the character
table information (see Theorem~\ref{thm:main}) and on the methods used to do
so; and to describe how such a check could be carried out independently by
anyone having available sufficient computing power.
\par

The \Atlas{} contains two essential pieces of information about each simple
group $S$ it lists: the ordinary character tables of all (or most) bicyclic
extensions $M.S.A$ of $S$, and all (or most) of the maximal subgroups of all
subgroups $S.A$ of the automorphism group of $S$ that contain $S$. We discuss
only these two pieces of data, and do not consider the other information also
given in the \Atlas, such as constructions of the groups and their
presentations, since we consider these two as the information most widely
used. As we will explain below, the situation for published proofs and
references for maximal subgroups is very satisfactory, so we mostly deal with
the question of (re)constructing the ordinary character tables. Here, we will
not try to follow the original proofs but rather give a modern approach,
which builds on the use of computer programs and is highly automatic (and so,
we claim, is much more reliable than hand calculations).
\par

One further comment is in order: nowadays, all the character tables contained
in the \Atlas, incorporating the corrections, and many more, are stored
electronically in the character table library \cite{CTblLib} of the computer
algebra system {\GAP} \cite{GAP}. Our checking will be with respect to these
electronic tables.
\par

The \Atlas{} tables in {\GAP}'s character table library have been constructed
from the data files (in the so-called ``Cambridge format'') which were used
also to create the printed \Atlas{}; hence the
ordering of rows and columns in the {\GAP} tables coincides with the \Atlas{}
ordering. We do not know how far the old ``Cambridge format'' files
represent exactly the contents of the printed \Atlas{}.
It might be possible to scan the printed \Atlas{} and to turn the
result into \GAP{} files; these could then be compared with the {\GAP} tables;
but we did not consider worthwhile following this cumbersome procedure here.
So for this practical reason all statements on accuracy and consistency will
only be relative to the electronic \GAP{} versions of the \Atlas, and we make
no claim on their agreement with the original printed version.
\par

The same data format has been used for the Brauer character tables in
\cite{JLPW95}: the information in \cite{JLPW95} depends on the \Atlas{}, the
Brauer character tables can be understood only relative to the ordinary
\Atlas{} character tables. When the Brauer character tables in \cite{JLPW95}
were prepared for both the printed version and their inclusion in {\GAP}'s
character table library, the ordinary \Atlas{} tables were already checked
systematically, and many of the errors listed in \cite{BN95} were found by
these checks.

\medskip
\noindent
{\bf Acknowledgement:} We thank Jean-Pierre Serre for raising the
question of reliability of \Atlas{} information, which led to the current
paper, and for comments on a preliminary version.

\section{Maximal subgroups}   \label{sec:maximal}

We begin by discussing the information in the \Atlas{} concerning maximal
subgroups of simple groups and of their automorphism groups. The situation
here is rather favourable, in the sense that published proofs for all of the
lists of maximal subgroups as printed in the \Atlas{} (modulo the corrections
listed in \cite{AtlCorr}) are available. For the sporadic simple groups and
their automorphism groups, references for the stated results are already
given in the original version of the \Atlas, and new information obtained since
then is referenced in \cite{AtlCorr}.
See \cite{W16} for a survey of this topic.

\par
As for the alternating groups $\fA_n$, $n\le13$, the question clearly is about
their primitive (maximal) subgroups, and these are well-known, see
e.g.~\cite{Mi97} for a classical reference, or \cite{Sims}.
The simple groups of Lie type in the \Atlas{} are of small Lie rank, and most
are of classical type. For the latter, much information on maximal subgroups
had already been accumulated in pre-\Atlas{} time, starting with the work of
Galois; the recent book by Bray, Holt and Roney-Dougal \cite{BHR} gives
complete proofs of the lists of maximal subgroups, and does \emph{not} rely
on the \Atlas{} lists. For the series of exceptional groups of types $\tw2B_2$,
$\tw2G_2$, $G_2$ and $\tw3D_4$, there exist published proofs \cite{Kl1,Kl2,Sz};
for $\tw2F_4(2)'$ and $F_4(2)$ proofs can be found in \cite{NW89,W84}.
Finally, for $\tw2E_6(2)$ the \Atlas{} does not claim to give complete
information.
\par
Thus, complete independent proofs for the maximal subgroup information in the
\Atlas{} are now available in the literature.

\section{Character tables}   \label{sec:chartabs}

We now turn to the more problematic question of character tables and their
correctness.

The tables for alternating groups, symmetric groups and their covering groups
are known by published classical work of Frobenius and Schur (see \cite{Sch11}
and the references therein).
\par

For many of the sporadic groups the \Atlas{} references published proofs,
for example a paper of Frobenius \cite{Fro} for the Mathieu groups.
(As Serre pointed out to us, while Frobenius sketched how the
character tables were constructed, he said nothing about
conjugacy classes. This is interesting, because even the existence
of $M_{24}$ was not completely clear at that time.)
However, such references do not exist for all of the sporadic groups.
\par

For the groups of Lie type, the
situation is even more opaque. No references are given in the \Atlas. While
there are published tables for some series of groups (for example for
$\SL_2(q)$ by Schur \cite{Sch07}, for $\tw2B_2(q^2)$ by Suzuki \cite{Sz},
and for $\SL_3(q)$ and $\SU_3(q)$ by Simpson and Frame \cite{SF73}, to mention
just a few), most of the tables for the groups of larger rank and in particular
for their decorations were computed using {\it ad hoc} techniques by
the \Atlas{} authors. Even today, the representation theory of finite groups
of Lie type, despite the tremendous achievements of George Lusztig, is not
capable of predicting the complete character tables of all the groups listed
in the \Atlas, in fact not even the character degrees in all cases.
\par

\subsection{Known mistakes}\label{sec:chartabs-mistakes}
Let us start by discussing the nature of known mistakes:
there are 142 entries marked as errors (three stars) in \cite{BN95},
and 17 such entries in \cite{AtlCorr}.
Many of the errors affect the descriptions of group constructions or maximal
subgroups, or indicator signs or power maps, and 27 concern character values.
Some of them could be fixed also by just changing power maps and some can be
detected by testing orthogonality relations --- for example, five columns of
the character table of $6.Fi_{22}.2$ were not printed, and a character value
$-1$ of $G_2(4)$ was listed as $1$. Some errors concern inconsistencies among
several characters. Consider for example the error on page 128 (the group in
question is the sporadic group $Suz$): the improvements list states
``Change sign of $i3$ in $\chi_{7}$, $\chi_8$, $\chi_{18}$, $\chi_{19}$,
$\chi_{21}$, $\chi_{22}$ on $6B$, $6C$''. For the simple group, one could
keep the character values, and adjust the power maps instead. However,
then one would have to change character values in central extensions of $Suz$.
For $G = 3.\PSU_3(8).3_1$ and $G = 3.\PSU_3(8).3_2$,
the problem was the consistent choice of irrationalities in the
faithful characters on the outer classes --- extensions to
$G$ of some faithful irreducible characters of the derived subgroup were
multiplied by $9$-th roots of unity, thus each of the shown characters exists
for a suitable group of the given structure but they do not fit to the same
isomorphism type.
But there are tables of (non-simple) groups which are wrong in a more serious
way, in the sense that characters were listed that cannot exist.
A wrong splitting of classes is shown for $2.\PSL_4(3).2_3$;
for both $G=\PSU_3(11).2$ and $G=2.\PSU_6(2).2$, the
extensions to $G$ of two irreducible characters of different degree of the
derived subgroup were interchanged. These mistakes are consistent
with the orthogonality relations and so are much harder to spot.

\subsection{Recomputing tables automatically}
We now propose our approach to reconstructing most of the character tables in
the (electronic version of the) \Atlas{} in a reproducible and essentially
automatic way. It relies on a powerful algorithm by Unger \cite{U06}.
We use his implementation which is available in \Magma \cite{Magma};
it uses no precomputed tables and does not rely on \Atlas{} bounds or data.
In the current version of {\GAP}~\cite{GAP},
the standard method to compute a character table is the less powerful
Dixon--Schneider algorithm \cite{Sch}. Both take as input a faithful
representation of a finite group, either as a permutation group or as a
matrix group over some finite field, and automatically compute the ordinary
character table of that group, including in particular the list of conjugacy
classes, their sizes and the power map on the classes.
\par

Now assume that we want to reconstruct the character table information for
a finite simple group $S$ appearing in the \Atlas. We proceed as follows.
First, the size and structure of the automorphism group $\Aut(S)$, the Schur
multiplier $M(S)$, and the action of the first on the second are well-known;
see, for example, \cite{GLS}. From this it is possible to compile
a list of all \emph{bicyclic} extensions $G=M.S.A$ for $S$ as considered in
the \Atlas: namely, both $M\le M(S)$ and $A\le\Out(S)$ are cyclic and $G$ is an
extension by $A$ of the central extension $M.S$ of $S$.

\subsection{Enumerating bicyclic extension}
Let $G$ be a group that contains normal subgroups $M < N$, and consider the
set of subquotients of the form $U / K$ with the property that $N \leq U$,
$K \leq M$, $K$ is normal in $U$, and both $U / N$ and $M / K$ are cyclic.
The group $G / N$ acts on this set by conjugation, and a set of class
representatives under this action contains all bicyclic extensions of $N / M$
that occur as subquotients of $G$, up to isomorphism. (Some representatives
may in fact be isomorphic; if we are interested in representatives up to
isomorphism, we must check this case by case.)

If $N$ is a Schur cover of a finite simple group $S$, so $M$ is
the Schur multiplier of $S$, and $G / M$ is isomorphic to the automorphism
group of $S$, then a set of class representatives yields all bicyclic
extensions of $S$, up to isoclinism. (Again, the set may be too large.)
We discuss the three most complicated cases occurring in the \Atlas{} in more
detail.

\subsubsection{Case 1: $S = \PSL_3(4)$}

The Schur multiplier $M$ and the outer automorphism group $A$ of $S$ have the
structures $3 \times 4^2$ and $D_{12}$ (the dihedral group of order twelve),
respectively. A group $G$ of the structure $M.S.A$ as mentioned above exists.
Since the subgroups $M_1$ and $M_2$ of order three and $16$ in $M$
are characteristic in $G$, we may consider the bicyclic extensions of $S$
that occur as subquotients of $G / M_1$ and $G / M_2$, and then get the
general bicyclic extensions of $S$ that occur as subquotients of $M$ as
subdirect products.

First we fix the notation for the cyclic subgroups of $G / N$. The unique
cyclic subgroup of order six is called $6$ by the {\Atlas}, its subgroups of
order two (the centre of the dihedral group) and three are called $2_1$
and $3$, respectively, and representatives of the other conjugacy classes of
subgroups of order two are called $2_2$ and $2_3$.

\begin{itemize}
\item
  The group $G / M_1$ has the structure $4^2.S.D_{12}$. Let $M / M_1$ be
  generated by commuting elements $a$, $b$ of order four, and let
  $c = (a b)^{-1}$. As stated in~\cite[p.~23]{Atl}, the outer
  automorphism group $G / N$ of $S$ acts as follows on $M / M_1$:
  \[
     \begin{array}{llll}
        6:    & a \mapsto b^3, & b \mapsto c^3, & c \mapsto a^3 \\
        2_2:  & a \mapsto   a, & b \mapsto   c, & c \mapsto   b \\
        2_3:  & a \mapsto a^3, & b \mapsto c^3, & c \mapsto b^3
     \end{array}
  \]
  The three subgroups of index two in $M / M_1$ are $\langle a, b^2 \rangle$,
  $\langle b, c^2 \rangle$, and $\langle c, a^2 \rangle$. Their normalisers
  in $G / N$ are the three Sylow $2$-subgroups. One of them contains the
  involutions $2_1$, $2_2$, $2_3$, thus we get the bicyclic extensions
  $2.S.2_1$, $2.S.2_2$, and $2.S.2_3$.

  (The other two Sylow $2$-subgroups of $G / N$ contain $2_1$ and conjugates
  of $2_2$ and $2_3$. Thus we get conjugate bicyclic extensions $2'.S.2_1$,
  $2'.S.2_2'$, $2'.S.2_3'$, $2''.S.2_1$, $2''.S.2_2''$, and $2''.S.2_3''$.)

  The group $G / N$ has the two orbits
  \[
    \left\{ \langle a \rangle, \langle b \rangle, \langle c \rangle \right\},
    \left\{ \langle a b^2 \rangle, \langle b c^2 \rangle,
              \langle c a^2 \rangle \right\}
  \]
  on the six cyclic subgroups of order four in $M / M_1$. We get two
  nonisomorphic central extensions of $S$ by a cyclic group of order four.
  Both extensions are normalised but not centralised by $2_1$, which inverts
  all elements in $M / M_1$.

  The second one, $(N / M_1) / \langle a b^2 \rangle$, is called $4_1.S$ by
  the {\Atlas}; it is centralised by $2_3$, its normaliser in $G/N$ is the
  elementary abelian group generated by $2_1$ and $2_3$. This yields the
  extensions $4_1.S.2_1$, $4_1.S.2_2$, and $4_1.S.2_3$.

  The first one, $(N / M_1) / \langle a \rangle$, is called $4_2.S$ by the
  {\Atlas}; it is centralised by $2_2$, its normaliser is the elementary
  abelian group generated by $2_1$ and $2_2$. This yields the extensions
  $4_2.S.2_1$, $4_2.S.2_2$, and $4_2.S.2_3$.

  (In both cases, the other two orbit points are stabilised by the other two
  Sylow $2$-subgroups of $G / M_1$, which yields the conjugate bicyclic
  extensions
  $4_1'.S.2_1$, $4_1'.S.2_2'$, $4_1'.S.2_3'$, $4_1''.S.2_1$,
  $4_1''.S.2_2''$, $4_1''.S.2_3''$,
  $4_2'.S.2_1$, $4_2'.S.2_2'$, $4_2'.S.2_3'$,
  $4_2''.S.2_1$, $4_2''.S.2_2''$, and $4_2''.S.2_3''$.)

\item
  The group $G / M_2$ has the structure $3.S.D_{12}$; the centraliser of
  $M / M_2$ in $G / N$ is the cyclic subgroup of order six,
  conjugation with the other elements of $G / N$ inverts $M / M_2$.

  Thus we get the bicyclic extensions $3.S$, $3.S.2_1$, $3.S.3$, $3.S.6$,
  $3.S.2_2$, and $3.S.2_3$, and their factor groups
  $S$, $S.2_1$, $S.3$, $S.6$, $S.2_2$, $S.2_3$,
  which are pairwise nonisomorphic.

\item
  Putting the pieces together, we get also bicyclic extensions in which the
  cyclic normal subgroup has order $6$ or $12$. Each of the above extensions
  with normal cyclic subgroup of order two or four and commutator factor
  group acting like one of the seven involutions can be combined with an
  extension with normal cyclic subgroup of order three and the same action
  of the commutator factor group.

  In summary, we get the following pairwise nonisomorphic bicyclic extensions:
  $6.S$, $6.S.2_1$, $6.S.2_2$, $6.S.2_3$,
  $12_1.S$, $12_1.S.2_1$, $12_1.S.2_2$, $12_1.S.2_3$,
  $12_2.S$, $12_2.S.2_1$, $12_2.S.2_2$, $12_2.S.2_3$.
\end{itemize}

\subsubsection{Case 2:  $S = \PSU_4(3)$}

The Schur multiplier $M$ and the outer automorphism group $A$ of $S$ have the
structures $3^2 \times 4$ and $D_8$ (the dihedral group of order eight),
respectively. A group $G$ of the structure $M.S.A$ as mentioned above exists.
Since the subgroups $M_1$ and $M_2$ of order four and nine in $M$ are
characteristic in $G$, we may consider the bicyclic extensions of $S$ that
occur as subquotients of $G / M_1$ and $G / M_2$, and then get the general
bicyclic extensions of $S$ that occur as subquotients of $M$ as subdirect
products.

First we fix the notation for the cyclic subgroups of $G / N$.
The unique cyclic subgroup of order four is called $4$ by the {\Atlas}, its
subgroup of order two (the centre of the dihedral group) is called $2_1$, and
representatives of the other conjugacy classes of involutions are called
$2_2$ and $2_3$.

\begin{itemize}
\item
  The group $G / M_1$ has the structure $3^2.S.D_8$.
  We identify $M / M_1$ with a $2$-dimensional vector space over $\FF_3$.
  The action of $G / N$ on this vector space is given by the matrices
  \[
     \left[ \begin{array}{rr} -1 & 0 \\ 0 & 1 \end{array} \right],
     \left[ \begin{array}{rr}  0 & 1 \\ 1 & 0 \end{array} \right]
  \]
  for the involutions from $2_2$ and $2_3$, respectively. The action of
  $2_1$ is given by the square of their product, which is the negative of
  the identity matrix; thus $2_1$ inverts all elements in $M / M_1$.
  The group $G / N$ has the two orbits
  \[
      \left\{ \pm [ 1, 0 ], \pm [ 0, 1 ] \right\},
      \left\{ \pm [ 1, 1 ], \pm [ 1, -1 ] \right\}
  \]
  on the nonidentity elements of $M / M_1$. We get two nonisomorphic central
  extensions of $S$ by a cyclic group of order three.

  The first one, $(N / M_1) / \langle [ 1, 0 ] \rangle$, is called $3_1.S$ by
  the {\Atlas}; it is centralised by $2_2$ and normalised by the elementary
  abelian group generated by $2_1$ and $2_2$. The third subgroup of order two
  in this subgroup is called $2_2'$, it centralises the conjugate extension
  $(N / M_1) / \langle [ 0, 1 ] \rangle$, which is called $3_1'.S$.

  The second one, $(N / M_1) / \langle [ 1, 1 ] \rangle$, is called $3_2.S$
  by the {\Atlas}; it is centralised by $2_3$ and normalised by the elementary
  abelian group generated by $2_1$ and $2_3$. The third subgroup of order two
  in this subgroup is called $2_3'$, it centralises the conjugate extension
  $(N / M_1) / \langle [ 1, -1 ] \rangle$, which is called $3_2'.S$.

  Thus we get the following pairwise nonisomorphic bicyclic extensions:
  $3_1.S$, $3_1.S.2_1$, $3_1.S.2_2$, $3_1.S.2_2'$,
  $3_2.S$, $3_2.S.2_1$, $3_2.S.2_3$ and $3_2.S.2_3'$.

  (Note that the centre of the groups $3_1.S.2_2$ and $3_2.S.2_3$ has order
  three, the other four groups have trivial centre. Analogously, the conjugate
  bicyclic extensions $3_1'.S.2_2'$ and $3_2'.S.2_3'$ have centres of order
  three, and the centres of $3_1'.S.2_1$, $3_1'.S.2_2$, $3_2'.S.2_1$, and
  $3_2'.S.2_3$ are trivial.)

\item
  The group $G / M_2$ has the structure $4.S.D_8$;
  the centraliser of $M / M_2$ in $G / N$ is the cyclic subgroup of order four,
  conjugation with the other elements of $G / N$ inverts $M / M_2$.

  Thus we get the bicyclic extensions $4.S$, $4.S.2_1$, $4.S.4$, $4.S.2_2$,
  and $4.S.2_3$,
  and their factor groups $2.S$, $2.S.2_1$, $2.S.4$, $2.S.2_2$, $2.S.2_3$,
  $S$, $S.2_1$, $S.4$, $S.2_2$, $S.2_3$, which are pairwise nonisomorphic.

\item
  Putting the pieces together, we get also bicyclic extensions in which the
  cyclic normal subgroup has order $6$ or $12$. Each of the above extensions
  with normal cyclic subgroup of order three and commutator factor group
  acting like one of the five involutions can be combined with an extension
  with normal cyclic subgroup of order two or four and the same action of the
  commutator factor group.

  In summary, we get the following pairwise nonisomorphic bicyclic extensions:
  $6_1.S$, $6_1.S.2_1$, $6_1.S.2_2$, $6_1.S.2_2'$,
  $6_2.S$, $6_2.S.2_1$, $6_2.S.2_3$, $6_2.S.2_3'$,
  $12_1.S$, $12_1.S.2_1$, $12_1.S.2_2$, $12_1.S.2_2'$,
  $12_2.S$, $12_2.S.2_1$, $12_2.S.2_3$, $12_2.S.2_3'$.

\end{itemize}

  The tables of $12_1.S.2_2'$ and $12_2.S.2_3'$ will be available in
  the next public release of \cite{CTblLib}.

\subsubsection{Case 3:  $S = \PSU_3(8)$}

\begin{itemize}
\item
  The outer automorphism group $A$ of $S$ has the structure $3 \times S_3$.
  The {\Atlas} notation for the four subgroups of order three in $A$ is
  $3_1$ for the central one, $3_2$ for the noncentral normal one,
  and $3_3$ and $3_3'$ for the remaining two, which are conjugate in $A$.
  In addition, $A$ contains nontrivial cyclic
  subgroups of order two and six, in each case three conjugates called
  $2$, $2'$, $2''$ and $6$, $6'$, $6''$, respectively.
  This yields the automorphic extensions $S$, $S.2$, $S.3_1$, $S.3_2$, $S.3_3$,
  and $S.6$.

\item
  The Schur multiplier $M$ of $S$ has order three,
  and a group of the structure $M.S.S_3$ exists.
  (There is no group of the structure $M.S.A$.)
  This yields the bicyclic extensions $3.S$, $3.S.3_1$, and $3.S.2$;
  the latter is conjugate and thus isomorphic to $3.S.2'$ and $3.S.2''$.
  Also a group of the structure $3.S.3_2$ exists.
  Concerning groups of the structure $3.S.6$,
  there is one such group that contains a given $3.S.3_1$ type group as
  a subgroup of index two.
  As stated in~\cite{BN95}, the groups $3.S.6'$ and $3.S.6''$, whose
  existence is claimed in~\cite[p.~66]{Atl}, are the extensions of the
  isoclinic variants of $3.S.3_1$ by their unique outer automorphisms
  of order two.
  Thus $3.S.6'$ and $3.S.6''$ are not isomorphic and in particular
  not conjugate to $3.S.6$ in some larger group.
  We can, however, ignore them since we are interested in the bicyclic
  extensions only up to isoclinism.
\end{itemize}

\subsection{Our algorithm}
For each bicyclic extension $G$ of a simple \Atlas{} group $S$, we
proceed as follows:
\begin{enumerate}
\item[(1)] We construct, or find in an existing database, a faithful
representation $H$ of $G$. For groups of Lie type, it is often easy to
construct such a representation from its very definition; both \GAP{} and
\Magma{} provide access to natural representations. At this point, we do not
actually need to establish that $H\cong G$. This will only occur at a later
stage.
\item[(2)] We then give this representation to the character table algorithm
in the computer algebra system of our choice (in our case \Magma{}).
\item[(3)] From the output (the character table of $H$) we can read off the
composition factors of $H$, assuming the classification of finite simple
groups, its derived subgroup and the centre of that subgroup.
The very few cases of simple
groups of identical order can easily be distinguished by additional information
on centraliser orders, for example. In many cases, this will already prove
that $H$ is isomorphic to $G$.
\item[(4)] It is then an easy computer algebra problem to check whether this
newly computed table is permutation isomorphic to the stored \Atlas{} table
for a group with this name, including the stored power map.
\end{enumerate}

We were able to apply this strategy to all but four simple groups $S$
contained in the \Atlas{} and found no discrepancy with the stored tables.
We summarise our main result.

\begin{theorem}   \label{thm:main}
 Let $G$ be a bicyclic extension of a simple group whose character table is
 given in the \Atlas, and different from $J_4$, $2.\tw2E_6(2)$,
 $2.\tw2E_6(2).2$, $B$, $2.B$ and $M$. The character table of $G$ has been
 automatically recomputed and found to agree with the one stored in the
 character table library of \GAP{}.
\end{theorem}

(See the remarks in the introduction about the relationship between the
\GAP{} tables and the printed \Atlas{} tables.)

Details of the computations underpinning this theorem can be found at
\cite{Br16b}, including the group generators that were used. The character
tables were computed using {\sc Magma} 2.21-4 on a 2.9GHz machine with 1TB Ram.
Most of the constructions were routine and used few resources, both in time
and memory; those with composition factors $\tw2E_6(2)$, $F_3$ and $Th$ were
challenging, the last taking 988\, 923 seconds.

\subsection{$J_4$, $2.\tw2E_6(2)$, and $2.\tw2E_6(2).2$}
We were able to compute the character tables of $J_4$, $2.\tw2E_6(2)$,
and $2.\tw2E_6(2).2$ using a combination of standard character theoretic
methods (induction from subgroups, LLL reduction, and the enumeration of
orthogonal embeddings); published information about conjugacy classes and
subgroups; and character tables verified in Theorem~\ref{thm:main} for
particular subgroups. These calculations are described in full in~ \cite{Br16}.
Specifically, the following information suffices for the computations.
\begin{itemize}
\item[$J_4$:]
  The description of conjugacy classes of $J_4$ as given in~\cite{Jan76};
  the character table of the subgroup of type $2^{11}\!:\! M_{24}$;
  this subgroup is the unique primitive permutation group on $2^{11}$ points
  of that order that can be embedded into $J_4$.
\item[$2.\tw2E_6(2)$:]
  The outer automorphism group of $\tw2E_6(2)$, a symmetric group
  on three points, acts faithfully on the classes of $\tw2E_6(2)$;
  and $2.\tw2E_6(2)$ contains subgroups of type $2.F_4(2)$;
  and $\tw2E_6(2)$ contains subgroups of types $F_4(2)$,
  $Fi_{22}$, $3\!\times\!\PSU_6(2)$, and $\OO_{10}^-(2)$;
  the character tables of the preimages of these subgroups under the
  natural epimorphism from $2.\tw2E_6(2)$.
\item[$2.\tw2E_6(2).2$:]
  The character tables of subgroups of types
  $2\times F_4(2) \times 2$ and \mbox{$3\times 2.\PSU_6(2).2$}.
\end{itemize}
The three character tables agree with the corresponding tables in
{\GAP}'s character table library.

\subsection{$B$, $2.B$ and $M$}
The tables of the remaining three groups excluded in the theorem are out of
reach to our methods; their verification will be considered by others.

\begin{remark}
Note that the Frobenius--Schur indicators of characters are not stored in the
\GAP{} tables, but recomputed when needed. But the indicators contained in
the old ``Cambridge format'' files were checked at the time of their
conversion into \GAP{} and so no undocumented errors should exist.
\end{remark}

\begin{remark}
The \Atlas{} contains all bicyclic extensions of simple groups of Lie type
possessing an exceptional Schur multiplier, with the sole exception of some
extensions of $\tw2E_6(2)$. More precisely, none of the bicyclic extensions
with one of the extensions being of degree~3 are given. For many applications,
it is useful to know those character tables as well.
\par
To our knowledge, the current status for $S=\tw2E_6(2)$ is as follows.
Frank L\"ubeck has computed the character table of $3.S$ using character
theoretic methods: Deligne--Lusztig theory contributes some information about
faithful characters, and this suffices for completing the whole character
table. The table of $6.S$ can be computed from the tables of
$2.S$ and $3.S$ automatically; the usual heuristics --- form tensor products
and apply LLL reduction --- is surprisingly successful.
Computing the tables of $3.S.2$ and $6.S.2$ is even easier since the outer
automorphism acts nontrivially on the centre. The tables of $2.S$, $6.S$,
$3.S.2$, and $6.S.2$ are available in~\cite{CTblLib}.
\end{remark}


\end{document}